\documentstyle[12pt]{article}

\begin{document}

\newtheorem{theorem}{Theorem}[section]
\newtheorem{corollary}{Corrolary}[section]
\newtheorem{lemma}{Lemma}[section]
\newtheorem{conjecture}{Conjecture}[section]
\newtheorem{convention}{Convention}[section]
\newtheorem{remark}{Remark}[section]
\newtheorem{definition}{Definition}[section]
\newtheorem{proposition}{Proposition}[section]
\renewcommand{\baselinestretch}{1.4}
\newcommand{\binom}{ { { {\rm {wt}}  a} \choose i } }
\newcommand{\wta}{{\rm {wt} }  a }
\newcommand{\R}{\mbox {Res} _{z} }
\newcommand{\wtb}{{\rm {wt} }  b }
\newcommand{\bea}{\begin{eqnarray}}
\newcommand{\eea}{\end{eqnarray}}
\newcommand{\be}{\begin {equation}}
\newcommand{\ee}{\end{equation}}
\newcommand{\g}{\bf g}
\newcommand{\hg}{\hat {\bf g} }
\newcommand{\h}{\bf h}
\newcommand{\hh}{\hat {\bf h} }
\newcommand{\n}{\bf n}
\newcommand{\hn} {\hat {\bf n}}
\baselineskip=12pt

\title{ {Vertex operator algebras associated to
modular invariant representations for $A_1 ^{(1)}$ }
\author{Dra\v{z}en Adamovi\'{c} and Antun Milas}}
%\pagestyle{myheadings}
%\markboth {Adamovi\'{c}, Milas}
%{Vertex operator algebras for $A_1 ^{(1)}$}
\pagenumbering{arabic}
\maketitle

\begin{abstract}
\baselineskip=12pt
We investigate vertex operator algebras $L(k,0)$ associated
with modular-invariant
representations for an affine Lie algebra $A_1 ^{(1)}$ ,
where k is 'admissible' rational number.
We show that VOA  $L(k,0)$ is rational in the category $\cal O$
and find all irreducible representations in  the
category of weight modules.
\end{abstract}
\section{Introduction}

Vertex operator algebras (VOA) are mathematical countrapart of conformal
field theory (CFT). It is very interesting that some representations of
affine Lie algebras carries the structure of VOA (or modules for VOA)
[FLM], [FZ], [MP].

The new sight in the theory of representation of VOA was made by Frenkel
and Zhu (see [FZ], [Z]) by introducing the associative algebra $A(V)$
associated
to VOA $V$. So called $A(V)$-theory gave an (theoreticaly) elegant way for the
classification of all ireducible representations of $V$ and for calculating
the 'fusion rules'. They also introduce the term of rational VOA which
is a VOA with finite number of irreducible modules, such that
every finitely generated module is completly reducible.

In [FZ], [MP], [KWn] the irreducible representations of the VOA
$L(k,0)$ ,$k \in {\bf N}$, associated to the irreducible highest weight
representations for an affine Lie algebra,
 were classifed. It seems that this case is much simplier because the
associative algebra $A(L(k,0))$ is finite dimensional (see [KWn]).

The main goal of this paper is a classification of the
irreducible representations
of the simple vertex operator algebra $L(k,0)$ for
$A_1^{(1)}$
on  the admissible rational level $k$.
Our main result is that irreducible $L(k,0)$--modules from the category
${\cal O}$ are exactly modular invariant representations for $A_1^{(1)}$.
To show this, we use $A(V)$--theory and identify $A(L(k,0))$ with the
certain quotient  of $U({\g}) $. Here we used  Malikov-Feigin-Fuchs formula
for the singular vectors in the Verma modules.
Then, by using classification of the irreducible representations in the
category ${\cal O}$,
we find all irreducible representations in the category of weight modules
for $A_1^ {(1)}$.

Feigin and Malikov
in [FM] have the geometrical approach to the similar problem
(see also [AY]). They
calculated conformal blocks for three admisseble modules associated to
 three different points on $\bf{CP^{1}}$. We interpret our result of
the classification of irreducible modules  in the terms of conformal blocks
considered in [FM].

\section{Preliminaries}

\subsection{Vertex operator algebras and modules }

\begin{definition} \label{def.v}
A {\em vertex operator algebra} is a
${\bf Z}$--graded vector space $V=\bigoplus_
{n \in{\bf Z}}V_n$ with
a sequence of linear operators
$\{ a(n) \mid n\in  {\bf Z} \}\subset End\ V$ associated to every $a\in V$,
such that for fixed
$a,b\in V$, $a(n)b=0$ for $n$ sufficently large.
We call the
 generating series
$Y(a,z)=\sum_{n\in {\bf Z}}a(n)z^{-n-1}\in (End\ V)[[z,z^{-1}]]$,
 {\em vertex operators} associated to $a$,
satisfy the following axioms:
\item[{\bf (V1)}]  $Y(a,z)=0$ iff $a=0$.

\item[{\bf (V2)}] There is a {\em vacuum} vector, which we denote by $\bf1$,
such that
$$Y({\bf 1},z)=I_V\,(I_V\, \mbox{is the identity of}\,\, End\,V). $$

\item[{\bf (V3)}]
There is a special element $\omega \in V$
(called the  {\em Virasoro\ element}), whose
 vertex operator we write in the form
$$Y(\omega,z)=\sum_{n\in {\bf Z}}\omega(n)z^{-n-1}=
\sum_ {n\in {\bf Z}}L_n z^{-n-2},$$ such that
$$L_0\mid_{V_n}=nI\mid_{V_n},$$
\bea
Y(L_{-1}a,z)=\frac{d}{dz}Y(a,z)\ \ \mbox{for every } a \in V ,
\eea
\bea
[L_m,L_n]=(m-n)L_{m+n}+\delta_{m+n,0}\frac{m^3-m}{12}c, \label{eq_vir}
\eea
where $c$ is some constant in ${\bf{C}}$,
which is called the {\em rank} of $V$.

\item[{\bf (V4)}] The {\em Jacobi identity} holds, i.e.
\bea
& &z_{0}^{-1}\delta\left(\frac{z_{1}-z_{2}}{z_{0}}\right)Y(a,z_{1})Y(b,z_{2})
-z_{0}^{-1}\delta\left(\frac{-z_{2}+z_{1}}{z_{0}}\right)Y(b,z_{2})Y(a,z_{1})
\nonumber\\
&=&z_{2}^{-1}\delta\left(\frac{z_{1}-z_{0}}{z_{2}}\right)Y(Y(a,z_{0})b,z_{2})
\eea
for any $a, b \in V$.
\end{definition}

The subspace $I$ of $V$ is called ideal if
$
Y(a,z)b \in I[[z,z^{-1}]]$
for every $a \in V, b \in I$.
Given an ideal $I$ in $V$ such that $\bf1 \notin I$, $\omega \notin I$,
the quotient $V/I$ admits a natural VOA structure( see [FZ] ).

\begin{definition}
Given an $VOA$  $V$, {\em a representation of $V$
(or V-module)} is a ${\bf Z}_{+} $-graded vector space
 $M=\bigoplus_{n \in {\bf Z}_{+}}M_{n}$ and
a linear map $$V\longrightarrow (End\ M)[[z,z^{-1}]],$$
$$a\longmapsto Y_M (a,z)=\sum_{n\in {\bf Z}} a(n)z^{-n-1},$$
satisfying
\item[{\bf (M1)}]  $a(n)M_m \subset M_{m+\deg a -n-1}\,\,\, \mbox{for every
homogeneous element a.}$
\item [{\bf (M2)}]  $Y_M(1,z)=I_M,$ and setting
$Y_M(\omega,z)=\sum_{n\in {\bf Z} }L_nz^{-n-2},$ we have
\bea
&&[L_m,L_n]=(m-n)L_{m+n}+\delta_{m+n,0}\frac{m^3-m}{12}c,\nonumber \\
&&Y_M(L_{-1}a,z)=\frac{d}{dz}Y_M(a,z)\,\, \nonumber
\eea
for every $\,\,\,a \in V.$

\item[{\bf (M3)}]  The {\em Jacobi identity} holds, i.e.
\bea
&&z_{0}^{-1}\delta\left(\frac{z_{1}-z_{2}}{z_{0}}\right)
Y_M (a,z_{1})Y_M (b,z_{2})
-z_{0}^{-1}\delta\left(\frac{-z_{2}+z_{1}}{z_{0}}\right)
Y_M(b,z_{2})Y_M(a,z_{1})
\nonumber \\
&&z_{2}^{-1}\delta\left(\frac{z_{1}-z_{0}}{z_{2}}\right)
Y_M(Y(a,z_{0})b,z_{2})
\eea
for any $a, b \in V$.
\label{definition_module}
\end{definition}

The  submodules, quotient modules, irreducible modules, completely
reducible modules are defined in the usual way ([FHL]).

\subsection{ Associative algebra $A(V)$ }

 Let $V$ be a VOA. For any homogeneous
element $a\in V$ and for any $b\in V$, following [Z] we define
\begin{eqnarray}
a* b={\rm Res}_{z}\frac{(1+z)^{{\rm wt a}}}{z}Y(a,z)b.
\end{eqnarray}
Then extend this product bilinearly to the whole space $V$. Let $O(V)$
be the subspace of $V$ linearly spanned by the elements of type
\begin{eqnarray}
{\rm Res}_{z}\frac{(1+z)^{{\rm wt a}}}{z^{2}}Y(a,z)b\;\;\;\mbox{for
homogeneous elements }a,b\in V.
\end{eqnarray}
Set $A(V)=V/O(V)$. The multiplication $*$ induces the multiplication on
the $A(V)$ and  $A(V)$ becomes an associative algebra. The image of $\bf 1$
in $A(V)$ becomes the identity element till the image of $\omega$ is in center
of
$A(V)$ (see [Z]).
Let $M=\oplus _{n \in {\bf Z} _+ }  M_{n} $ be a $V$--module. For a homogeneous
element $a \in V$ we define $o(a)=a(\deg a-1)$. From the definition of $M$
follows that operator $o(a)$ preserves the graduation of $M$.
\begin{theorem} \label{fz}
\item [(a)]
On $\mbox {End}(M_0)$ we have
\bea
&&o(a)o(b) = o(a*b) \nonumber \\
&&o(x) =0 \nonumber
\eea
for every $a,b \in V$, $x \in O(V)$.  The top level $M_0$ is an
$A(V)$--module.
\item[(b)]  Let $U$ be an $A(V)$--module, there exists $V$--module
$M$ such that $A(V)$--module $M_0$ and $U$ are isomorphic.

Thus, we have one-to-one correspondence between irreducible $V$--modules
and irreducible $A(V)$--modules.
\end{theorem}

We have the following consequence of the definition of
$A(V)$.
\begin{proposition}  \label{id}
Let $I$ be an ideal of $V$
Assume ${\bf 1} \not \in I,
\omega \not \in I$. Then the associative algebra
$ A(V/I)$ is isomorphic to $A(V) / {[I]}$, where $[I]$ is the image
 of $I$ in $A(V)$.
\end{proposition}

\subsection {Vertex operator algebras associated to affine Lie algebras}
Let ${\g}$ be a finite-dimensional simple Lie algebra over ${\bf C}$.
The affine Lie algebra
${\hg}$ associated with ${\g}$ is defined as
$
{\g} \otimes {\bf C}[t,t^{-1}] \oplus {\bf C}c
$
with the usual commutation relations.
Let  ${\g} = {\n}_- + {\h} + {\n}_+$ and  ${\hg}={\hn}_- +{\hh} + {\hn}_+$
be the usual triangular decompositions for ${\g}$ and ${\hg}$ and
 $P = {\bf C} [t] \otimes {\g} \oplus {\bf C}c $ be upper parabolic
subalgebra.
Let $U$ be any ${\g}$--module. Considering $U$ as a $P$--module
, we have the induced module ( so called generalized Verma module)
$
M({\ell},U) = U({\hg} ) \otimes _{U(P)} U $, where the central element $c$ acts
as multiplication with ${\l} \in {\bf C}$.

For $\lambda \in {\h} ^*$ with $M(\lambda) $ we denote Verma module and
with $V(\lambda)$ its  irreducible quotient.

Set $M(\ell,\lambda)=M(\ell, V(\lambda))$. Let
$L(\l,\lambda)$ denotes its irreducible quotient.

\begin{theorem} ([FZ]) Every $M(\ell,0)$ $\ell \ne -g$ (where
$g$ denotes dual Coxeter number) has the structure of VOA.
Let $U$ be any ${\g}$--module. Then
every $M(\ell, U)$ is a module for $M(\ell,0)$. In particulary
$M(\ell, \lambda)$
is $M(\ell,0)$--module.
\end{theorem}

\begin{theorem} \label{nkl}
The associative algebra $A(M(\ell,0))$ is
canonically isomorphic to  $U(\g)$ and the isomorphism
$F  : A(M({\ell},0)) \rightarrow U({\g})$
is given by :
\bea \label{as.i}
  F  \ \ [a_1(-i_1 -1) \cdots a_n(-i_n -1) {\bf 1}]
  = (-1) ^{i_1 + \cdots i_n}  a_n \cdots a_1 .
\eea
for every $a_1,\cdots,a_n \in {\g}$ and every $i_1,\cdots,i_n \in {\bf Z}_+$.
\end{theorem}

\section{ Irreducible modules for VOA $L(k,0)$ in the category ${\cal O}$ }

\subsection{Modular invariant representations for $A_1^{(1)}$}

Let now ${\g}=sl(2,\bf C)$ with generators $e,f,h$ and relations:
$[h,f]=-2f$, $[h,e]=2e$, $[e,f]=h$.
Let $\Lambda_0, \ \Lambda_1$ denote the fundamental weights for ${\hg}$,
and $\omega$ the fundamental weight for ${\g}$.

\begin{definition} $k=p/q \in {\bf Q}$ is called admissible if
$q\in {\bf N}$, $p\in {\bf Z}$, $(p,q)=1$ and $2q+p-2 \ge 0$.
\end{definition}

 In [KW] V.Kac and M.Wakimoto define modular invariant representations.
 They also define weights which have admissible level
and satisfy some technical conditions (for definition see [KW]).
They call them admissible weight.

The following proposition describes the admissible weights
and  modular invariant representations on level $k$:
\begin{proposition}
Let $k=p/q \ \in Q$ be admissible. Set $t=k+2$.
Define :
$$
P^{k} = \{ (k - n + mt) \Lambda_0 + (n - mt ) \Lambda_1 ,
\ m,n \in {\bf Z}_+ , \ n \le 2q+p-2, \ m \le q-1 \}.
$$
Let $M$ be any irreducible highest weight module , with the highest weight
$\lambda$.
The following statments are equvivalent :
\item[(1)] $M$ is a modular-invariant.
\item[(2)] $\lambda$ is  an admissible  weight.
\item[(3)] $\lambda \in  P^{k} $.
\end{proposition}
(For proof  see [KW]).

We will need the following describtion of the set $P^{k}$.
\begin{lemma} \label{opis}
Let $\lambda \in {\hh} ^*$. Then $ \lambda \in  P^{k}$ if and only if
$$\langle \lambda, c \rangle =k, \
\ \langle \lambda , h \rangle =(N-it-j)$$
where  $i \in  \{ 0, \cdots ,l \}$ ,$ j \in \{ 1,\dots , N \}$,
 $N=2q+p-1, l=q-1. $
\end{lemma}
 By using Cor.2.1.in [KW] or Kac determinant formula we have

\begin{theorem}  \label{kacwak}
Let $k=p/q \in {\bf Q}$ be  admissible . Then
$$L(k,0)=M(k,0)/U ({\hg}) v_{sing},$$
where vector $v_{sing}$ is the unique singular vector of the weight
$k \Lambda_{0} - q(2q+p-1) \delta + (2q+p-1) \alpha$.
\end{theorem}
We also need the following theorem:
\begin{theorem} \label{kw}
(Kac--Wakimoto)
Let $M$ be a ${\hg}$-- module from the category $\cal O$ such that for any
irreducible subquotient $L(\mu )$ the weight  $\mu $ is admissible.
Then  ${\hg}$--modul $M$ is  completely reducible.
\end{theorem}

\subsection{Malikov-Feigin-Fuchs formula}

Now, we recall the result of Malikov-Feigin-Fuchs which give us the
singular vector in form with "rational powers"(see [MFF]).
\begin{theorem}
(Malikov-Feigin-Fuchs)
The maximal submodule of $M(k,0)$ is generated by the singular vector
given by the
$v_{sing} = F(k).{\bf 1}$ where
\be \label{s.form}
F(k) = e(-1)^{N+lt} f(0)^{N+(l-1)t} \cdots     f(0)^{N-(l-1)t}
 e(-1)^{N-lt},
\ee
for $N=2q+p-1$,$l=q-1$ and $t=p/q +2$.
\end{theorem}
\begin{remark}
In [MMF] were proved that this formula realy make sense, because only with
commuting we can transform formula (\ref{s.form}) in the usually form in
$U({\hg})$.
\end{remark}
\subsection{Fundamental lemma}

First we define  :
\bea
\epsilon : \ \ \
&&   U({\hn} _- ) \rightarrow U({\g}) \nonumber \\
&& a_1(-i_1)\cdots a_s(-i_s) \mapsto a_1\cdots a_s,  \nonumber
\eea
for every $a_1,\dots,a_s \in {\g}$, $s \in {\bf N} $.

In the same way as in [F] we have:
\begin{proposition} \label{Fuchs}
$$\epsilon(F(k) ) = \prod_{i=1}^{l} \prod_{j=1}^{N} p_{i,j}(h) e^{N},$$
where $p_{i,j}(h)=ef + (it+j-1)h - (it+j)(it+j-1)$.
\end{proposition}

We define the ${\bf Z}$--graduation on $U({\hg})$ :
\be
\deg a_1(-i_1) \cdots a_k(-i_k) = i_1 + \cdots + i_k ,
\ee
for every $   a_1, \dots , a_k \in {\g}$.

In the following lemma we will use the ordinary transposing $^T$ in $U({\g})$
(see [Dix]).

\begin{lemma} \label{kom.d}
Let $g \in U({\hn}_-)$, such that $\deg g = n$ . Then we have
$$
\epsilon(g) \equiv  (-1) ^ n  ( F [g.{\bf 1} ]) ^{T} \ \mbox{mod}
\  U({\g}) {\n}_- .
$$
\end{lemma}
{\em Proof.} First notice that ${\n} _{-} . {\bf 1}=0$.  Since $\deg g =n$ ,
one can write $g$ in a form
$$g=\sum_{i=0} ^{r} g_i f(0) ^ i,$$
where
$$g_i = \sum a_{i_1} ^{(i)} (-j_1-1) \cdots a_{i_t} ^{(i)} (-j_t-1),$$
$a_{i_1} ^{(i)},\dots , a_{i_t}  ^{(i)}\in {\g}$, $j_1,\dots ,j_t \in {\bf Z}
_+$,
$j_1+\cdots +j_t + t = n$, $r \in {\bf Z} _+$, and get
$$ \epsilon (g) \equiv g_0 \ \mbox{mod} \  U({\g}) {\n}_-. $$
Set $a_{i_j}= a_{i_j} ^{(0)}$.
Since
$$
g.{\bf 1} = g_0 . {\bf 1} =   \sum a_{i_1} (-j_1-1) \cdots a_{i_t}(-j_t-1),
$$
we have that :
\bea
F([ g.{\bf 1}]) ^ {T} =&&\sum (-1) ^{n-t} ( a_{i_t} \cdots a_{i_1} ) ^ T
\nonumber \\
= && (-1) ^n \sum  (a_{i_1} \cdots a_{i_t} ) \nonumber
\eea
and lemma holds.  ${ \ \ \Box}$

Set $Q=F([v_{sing}]) \in U({\g})$. From proposition \ref{Fuchs} and lemma
\ref{kom.d} we have
\begin{lemma}
$$Q^T \equiv (-1) ^ {q(2q-p-1)}
\prod_{i=1}^{l} \prod_{j=1}^{N} p_{i,j}(h) e^{N} \ \mbox{mod} \ U({\g}){\n}_-
$$
where polinomials $p_{i,j}$ are as in proposition \ref{Fuchs} .
\end{lemma}

\subsection{Classification of representation}

The vertex operator algebra $M(k,0)$ has the maximal ideal $M^1 (k,0)$.
It is generated by the vector $v_{sing}$. Let
 $L(k,0)$ be the quotient  VOA. The proposition \ref{id} and
theorem \ref{nkl} imply:

\begin{proposition}
$A(L(k,0))$ is isomorphic to
$U({\g})/I$ where $I$ is a two sided ideal generated by the vector $Q$.
\end{proposition}

Let $U$ be any  $A(L(k,0))$--module. Then $U$ is a ${\g}$--module.
We have
\begin{proposition} \label{keyprop}
Let $U$ be any $U({\g})$--module.
Then the following statments are equivalent:
\item[(1)]
$U$ is a $A(L(k,0))$--module,
\item[(2)]
$Q.U=0$.
\end{proposition}

 Set $R=U({\g}).Q$ and $R^T=U({\g}).Q^T$. Clearly
$R$ and $R^T$ are irreducible ${\g}$--modules and
$R \cong R^T \cong V(2 N\omega) \cong V^*(2 N\omega)$.

{}From this facts and proposition \ref{keyprop} one can  obtain:
\begin{lemma} Let $V(\mu)$ be the irreducible highest weight ${\g}$--module
with the highest weight vector $v_{\mu}$.
The following statments are equivalent :
\item[(i)] $V(\mu)$ is a $A(L(k,0))$--module ,
\item[(ii)] $RV(\mu)=R^T V(\mu)^*=0$,
\item[(iii)] $R_0 v_{\mu}=R_0^T v_{\mu}^*=0$,

where $R_0 \ (\ R_0^T)$ denotes the zero-weight subspace of $R\ (R^T) $.
\end{lemma}
For
$p \in  S({\h})$   and ${\mu} \in {\h} ^*$
 define $p ({\mu}) \in {\bf C}$ with
$p(h).v_{\mu}=p({\mu})v_{\mu}$.

Let $u_1 \in R_0$ and $u_2 \in R_0 ^T$.
Clearly there exists unique polynomials $p_1,p_2 \in S({\h})$
such that
$$u_1 \equiv p_1(h) \ \mbox{mod} \ U({\g}){\n}_+ \ \ \
u_2 \equiv p_2(h) \ \mbox{mod} \ U({\g}){\n}_-. $$
Then $u_1.v_ {\mu} = p_1 ({\mu}) v_{\mu}$ and
$u_2.v ^* _{\mu}= p_2 (-{\mu}) v ^* _{\mu}$.

 We have  :
\begin{lemma} \label{key}
There is one-to-one correspondence between  each two of the following three
sets :
\item[(1)]
$\mu \in {\h}^*$ such that $V(\mu)$ is $A(L(k,0))$--module ,
\item[(2)]
$\mu \in {\h}^*$ such that $p_1(\mu)=0$,
\item[(3)]
$\mu \in {\h}^*$ such that $p_2(-\mu)=0$.
\end{lemma}

\subsection{The main theorem}

The following lemma is obtained by direct calculation :
\begin{lemma} \label{pomoc}
$$ [f^N, ef+(it+j-1)(h-(it+j))]=(-N-1+it+j)(h-it-j+N)f^N $$
\end{lemma}

\begin{proposition} \label{class}
All irreducible $A(L(k,0))$--modules from the category $\cal O$ are:
$V(r \omega ), r \in S$,  where
\be
S=\{ N-it-j  \  : \ i=0,...,l; \ j=1,...,N \}.
\ee
\end{proposition}
{\em Proof.}
Let $u \in R_0^T$. Then
$u=(ad \ f) ^{N} .Q^T \equiv \ f^N Q^T \ \mbox{mod} \ U({\g}){\n}_-.$
By using lemma \ref{pomoc} we have
$$u  \equiv c_1\prod_{i=1}^{l} \prod_{j=1}^N
\ q_{i,j}(h)f^Ne^N \ \mbox{mod} \ U({\g}){\n}_-,$$
where $q_{i,j}(h)=h-it-j+N$, $c_1 \in {\bf C}$ .
Since $f^N e^N \equiv c_2 h(h+1) \cdots (h+n-1) \
\mbox{mod} \ U({\g}){\n}_-$, for some $c_2 \in {\bf C}$,
we conclude that polynomial
$p_2$ from lemma \ref{key} is proportional to
$$  \prod_{i=0}^{l} \prod_{j=1}^N (h-it-j+N).$$
Now, proposition follows from lemma \ref{key}.
${\ \ \Box}$

We can obtain the main theorem:
\begin{theorem} \label{main}
 The set $\{L(k,r\omega )\ : r \in S \}$ provides a complete
list of irreducible $L(k,0)$--modules from the category $\cal O$.
Moreover, the irreducible $L(k,0)$--modules from the category
 $\cal O$ are exactly irreducible highest weight representations with
 admissible highest weights.
\end{theorem}
{\em Proof.}  Proposition \ref{class}   and theorem \ref{fz} imply that
$L(k,r \omega )$, for $r \in S$ , are all irreducible $L(k,0)$--modules from
the
category ${\cal O}$. The secend statement follows from lemma \ref{opis}.
${\ \ \Box}$

\begin{theorem}
Let $M$ be a $L(k,0)$--module from the category $\cal O$. Then
$M$ is  completely   reducible $L(k,0)$--module.
\end{theorem}
{\em Proof.}
Let $M$ be a $L(k,0)$--module from the category $\cal O$ and let
$N$ be an irreducible subquotients of $M$. Then $N$ is irreducible
$L(k,0)$--module. From the theorem \ref{main} follows that
$N$ is a irreducible highest weight module with admissible highest weight.
Now theorem \ref{kw} implies that $M$ is completely reducible
${\hg}$--module and so completely reducible  $L(k,0)$--module.
${ \ \ \Box}$
\begin{remark}
 Vertex operator algebra is by definition
rational if it has only finitely many irreducible modules and if every finitely
generated module is a direct sum of irreducible ones.
We have showed that VOA $L(k,0)$ has finitely many irreducible modules in
the category  $\cal O$ and every module from the category $\cal O$  is
completely reducible.
 By using this arguments we  say  that
the vertex operator algebra $L(k,0)$, for $k \in {\bf Q}$ admissible, is
rational in the category $\cal O$.
\end{remark}
\begin{remark} In [A] were considered some modular invariant representations
for
$C_{\ell} ^{(1)}$  and proved that  the
VOA   $L(n-\frac 3 2, 0)$, $n \in {\bf N}$, is rational
in the category  $\cal O$.
\end{remark}
We have :
\begin{conjecture}
Let ${\g}$ be any simple finite-dimensional Lie algebra and $L(k,0)$
associated vertex operator algebra such that the highest weight of $L(k,0)$ is
admissible. Then $L(k,0)$ is rational in the category $\cal O$.
\end{conjecture}
\section{Irreducible modules for $L(k,0)$ in the category of
weight modules}

Let $M$ be any irreducible $L(k,0)$--module. From [FHL] we have that
the contragradient $L(k,0)$--module $M^*$ is also irreducible. Moreover,
 $M^{**}$ and $M$ are isomorphic $L(k,0)$--modules.
One can easily see that for $M=L(k,\lambda)$
$M^* _0$ is isomorphic to $V(\lambda) ^*$. We have  :
\begin{proposition}
$V(r \omega) ^* $,$r \in S$, are all irreducible lowest weight
$A(L(k,0))$-- modules.
\end{proposition}

Set $E_{r,\mu}=t^{\mu} \bf {C}[t,t^{-1}]$ where
$r,\mu \in{ \bf C}$
and $E_{i}=t^{\mu + i}$. We define $U({\g})$ action on $E_{r,\mu}$
by the following formulas :
\be \label{form}
e.E_{i}=-(\mu + i)E_{i-1},\ \ h.E_{i}=(-2 \mu -2i + r), \ \
f.E_{i}=(\mu +i-r)E_{i+1} .
\ee
We will find all  pairs $(r,\mu)$, such that $E_{r,\mu}$ is
irreducible $A(L(k,0))$--module.
\begin{theorem}
Set $T=
\{ \ (r,\mu) \ : \ r \in S - {\bf Z} _+ , \ \mu \notin \bf Z, \
r - \mu \notin \bf Z \}$.
Then $E_{r,\mu}$ is irreducible $A(L(k,0))$--module if and only if
$(r,\mu) \in T$.
\end{theorem}
{\em Proof.}
First, we notice that $E_{r,\mu}$ is an irreducible $U({\g})$--module
iff $\mu \notin \bf Z$ and $r - \mu \notin \bf Z$.

By using (\ref{form}) we have:
\be \label{form1}
Q.E_{i}=(p_{0}(r)+p_{1}(r)(i+\mu)+ \cdots +p_N(r)(i+\mu)
^N)E_{i-N},
\ee
for some polynomials $p_0,p_1 \cdots p_N$ and  $\mu \in \bf C$.

$Step \  1$. Let $E_{r,\mu}$ be an irreducible $A(L(k,0))$--module, then
$r  \in S - {\bf Z}_{+}$.

{}From the proposition \ref{keyprop} follows $Q.E_i=0$ for all $i \in \bf Z$.
{}From this fact and from (\ref{form1}) we have that:
$$ p_0(r)=p_1(r)= \cdots =p_N(r)=0 .$$
If $\mu=0$ we have that ${\bf C}[t]$ is a submodule of ${\bf C}[t,t^{-1}]$
isomorphic to $M(r\omega)$. From (\ref{form1}) and proposition \ref{keyprop}
follows that $M(r\omega)$ is $A(L(k,0))$--module.
Then the theorem \ref{main} implies that
$r  \in
S-\bf Z _+$ (in this case $V(r\omega)=M(r\omega)$).

$Step \ 2$. If $(r,\mu) \in T$ then $E_{r,\mu}$ is the irreducible
$A(L(k,0))$--module.

Since  $r \in S-{\bf Z} _+$, we have that
$M(r\omega)=V(r\omega)\cong {\bf C}[t]$ and $Q. {\bf C}[t] = 0$.
By using (\ref{form}) for $\mu =0$, we conclude   that
$p_0(r)=p_1(r)= \cdots=p_N(r)=0$. This fact implies that
$$Q.E_{i}=0 \ \ \ \mbox{for all}\ \ \  i \in \bf Z,$$
and we obtain that for $(r,\mu) \in T$,
$E_{r,\mu}$ is $A(L(k,0))$
--module. ${\ \ \ \Box}$

Recall that $U({\g})$--module $U$ is called weight module if ${\h}$ acts
semisimple on $U$ and all weight subspaces are finite-dimensional.
We know that irreducible weight modules
are : highest weight, lowest weight and modules $E_{r ,\mu}$ defined
by ( \ref{form} ) .

We have obtained :
\begin{corollary}
Let $U$ be an irreducible $A(L(k,0))$--weight module.
Then $U$ is one of the following modules :
\item[(1)]
$V(r\omega)$, $r \in S$,
\item[(2)]
$V(r\omega)^*$, $r \in S$ or
\item[(3)]
$E_{r,\mu}$, $(r,\mu) \in T$.
\end{corollary}

\begin{theorem}
Let $M$ be an irreducible $L(k,0)$-- module such that $M_0$ is a weight module.
Then $M$ is one of the following modules :
\item[(1)]
$L(k, V(r\omega))$, $r \in S$,
\item[(2)]
$L(k, V(r\omega)^*)$, $r \in S$ or
\item[(3)]
$L(k, E_{r,\mu})$, $(r,\mu) \in T$.
\end{theorem}

\section{ Connection with geometrical approach}

Let  $E$ be any finite subset of ${\bf CP^{1}}$ and let ${\bf g}(E)$ denotes
the
Lie algebra of meromorphic functions on  ${\bf CP^{1}}$ holomorphic outside
$E$ with values in ${\bf g}$. For every $z \in {\bf CP^{1}}$ and $\lambda \in
h^{*}$ we can define irreducible highest weight ${\bf g}(z)$--module
$L(k,\lambda,z)$ attached to $z$ (for definition see [FM]).

Let $z_{1}, z_{2}, z_{3}$ be three different points on ${\bf CP^{1}}$ and
$\lambda_{1},\lambda_{2},\lambda_{3} \in h^*$. We consider
${\bf g}(z_{1}, z_{2}, z_{3})$--module
$ L(k,\lambda_1,z_{1}) \otimes L(k,\lambda_2,z_{2}) \otimes
L(k,\lambda_3,z_{3})$
and the space of coinvariants
$$
H^{\circ}({\bf g}(z_{1}, z_{2}, z_{3}),\
L(k,\lambda_1,z_{1}) \otimes L(k,\lambda_2,z_{2})
\otimes L(k,\lambda_3,z_{3})).
$$

In previous section we showed that the irreducible $L(k,0)$--modules (in the
category ${\cal O}$) are exactly $L(k,r\omega)$, $r \in S$.
Those modules were considered in [FM].
They calculated dimension of the space
$$
H^{\circ}({\bf g}(0,1,\infty ),\
L(k,r_1\omega  ,0) \otimes L(k,r_2\omega,1)
\otimes L(k,r_3\omega,\infty))
$$
(this space is also called conformal block)
for all triples $r_{1},r_{2},r_{3} \in S$ and obtained "fusion algebra".

When $r_{3}=0$ their result implies  that
$$
\mbox{dim} \
H^{\circ}({\bf g}(0,1,\infty ),\
L(k,r_1\omega,0) \otimes L(k,r_2\omega,1)
\otimes L(k,0,\infty) )
=
\left \{
\begin{array}{ll}
1& \mbox{if} \ r_1=r_2  \in S \\
0&\mbox{otherwise}
\end{array} \right.
$$
We have the following characterisation of $L(k,0)$--modules:
\begin{theorem}
$L(k,\lambda)$ is  a $L(k,0)$--module if and only if
$$
\mbox{dim} \
H^{\circ}({\bf g}(0,1,\infty ),\
L(k,\lambda,0) \otimes L(k,\lambda,1)
\otimes L(k,0,\infty) ) = 1.$$
\end{theorem}

As in [FHL],  for three modules we can define so called fusion rules
(dimension of the space of intertwining operators).
{}From the previous theorem follows that when one of the modules is
$L(k,0)$ then
 fusion rules and   dimension  of corresponding conformal block are equal.
 It seems that this is thrue for any three modules.

Department of Mathematics, University of Zagreb,
Bijeni\v{c}ka 30, 41000 Zagreb, Croatia

e-mail : adamovic@cromath.math.hr, \  milas@cromath.math.hr

\end{document}